\documentclass[sn-mathphys]{sn-jnl}

\usepackage{lineno,hyperref,subfig,multirow,booktabs}
\usepackage{algorithm}
\usepackage{float}
\usepackage{amsmath}
\usepackage{xcolor}

\DeclareMathOperator*{\argmax}{argmax}

\jyear{2021}%

\theoremstyle{thmstyleone}%
\theoremstyle{thmstyletwo}%

\theoremstyle{thmstylethree}%

\begin{document}

\title[A Diversity-Aware Memetic Algorithm for the Linear Ordering Problem]{A Diversity-Aware Memetic Algorithm for the Linear Ordering Problem: Improving Best-Known Solutions for Standard Benchmarks}


\author[1]{\fnm{Lázaro} \sur{Lugo}}
\equalcont{These authors contributed equally to this work.}

\author*[1]{\fnm{Carlos} \sur{Segura}}\email{carlos.segura@cimat.mx}
\equalcont{These authors contributed equally to this work.}

\author[2]{\fnm{Gara} \sur{Miranda}}
\equalcont{These authors contributed equally to this work.}

\affil[1]{\orgdiv{Área de Computación}, \orgname{Centro de Investigación en Matemáticas (CIMAT)}, \orgaddress{\street{Callej\'on Jalisco s/n, Mineral de Valenciana}, \city{Guanajuato}, \postcode{36240}, \state{Guanajuato}, \country{México}}}

\affil[2]{\orgdiv{Departamento de Ingeniería Informática y de Sistemas}, \orgname{Universidad de La Laguna}, \orgaddress{\street{Avda. Astrofísico Francisco Sánchez, s/n}, \city{San Cristóbal de La Laguna}, \postcode{38200}, \state{Santa Cruz de Tenerife}, \country{España}}}


\abstract{The Linear Ordering Problem (LOP) is a very popular NP-hard combinatorial optimization problem with many practical applications that may require the use of large instances.
The Linear Ordering Library (LOLIB) gathers a set of standard benchmarks widely used in the validation of solvers for the LOP.
Among them, the xLOLIB2 collects some of the largest and most challenging instances in current literature. 
In this work, we present new best-known solutions for each of the 200 complex instances that comprises xLOLIB2.
Moreover, the proposal devised in this research is able to achieve all current best-known solutions in the rest of instances of LOLIB and improve them in other 93 cases out of 485, meaning that important advances in terms of quality and robustness are attained.
This important advance in the field of the LOP has been possible thanks to the
development of a novel Memetic Algorithm (MA) that was designed by taking into account some
of the weaknesses of state-of-the-art LOP solvers.
One of the keys to success is that the novel proposal allows for a gradual shift from exploration to exploitation, which is done by taking into account the stopping criterion and elapsed period of execution to alter the internal decisions taken by the optimizer. 
The novel diversity-aware proposal is called the Memetic Algorithm with Explicit Diversity Management (MA-EDM) and extensive comparisons against state-of-the-art techniques provide insights into the reasons for the superiority of MA-EDM.
}




\maketitle

\section{Introduction}\label{sec1}

The Linear Ordering Problem (LOP) is a classical combinatorial optimization problem that arose in the field of economics.
It has applications in other research areas, such as machine learning and logistics. 
The problem was initially introduced in 1936 by Leontief~\cite{Leontief:1936} and formalized in 1958 by Chenery and Watanabe~\cite{Chenery:58}, but it was not until decades later that Garey and Jonhson~\cite{Garey:79} demonstrated that LOP is an NP-hard problem.
Although the problem has multiple applications and formulations, 
it is quite common to define the problem as a matrix triangulation problem.

Given a matrix $M_{n \times n} = (m_{ij})$, the triangulation problem is to determine a simultaneous permutation $\sigma$ of the rows and columns of $M$ such that the sum of the entries above the main diagonal is maximized (or equivalently, the sum of subdiagonal entries is minimized)~\cite{Marti:2012}.
Hence, the LOP can be formulated as finding a permutation $\sigma$ that maximizes the following equation:

\begin{equation}
 \sum_{i = 1}^{n - 1} \sum_{j = i + 1}^{n} m_{\sigma_{i}\sigma_{j}}     
\end{equation} 

Several exact and approximate solvers have been proposed for this combinatorial optimization problem \cite{Marti:2011}. 
The interest in this problem has continued over the years, resulting in research compilations such as \cite{Marti:2012}. 
So far, different algorithms have been devised that yield very promising solutions for the LOP instances addressed in \cite{Marti:2012}. 
Within these proposals, the Memetic Algorithms (MAs) have proven to be very effective in the search for improved solutions to the LOP~\cite{Ceberio:2015}. 
As described in the following section, most trajectory- and population-based state-of-the-art approaches for the LOP rely on restarting procedures or other mechanisms to regain diversity or to abandon oversampled regions, in an effort to avoid stagnation and premature convergence.
In addition, they usually include problem-dependent information to exploit the particularities of the problem's structure.
The performance of schemes based on recovering diversity has been improved in several large-scale problems by methods that explicitly manage diversity~\cite{Crepinsek:13}, and methods that relate the amount of diversity to the elapsed period of execution and stopping criterion have been particularly successful~\cite{Segura2016}.
However, this last design principle has not been adopted in current LOP solvers.
The hypothesis of this research is that there is still a big room for improvement in the design of LOP solvers for large instances, and thus, the solutions that are used as a reference when validating LOP solvers might be enhanced significantly.
Knowing near-optimal values is important to properly measure the capabilities of different solvers.
Thus, this research aims to design a novel long-term solver that is able to obtain new best-known solutions for the most challenging instances of the LOP.
Given the promising behavior of MAs and the lack of population-based LOP solvers that establish a dependency between the amount of diversity maintained in the population and the elapsed execution period and stopping criterion, a novel MA that adopts this design principle is applied.

The novel diversity-aware memetic algorithm proposed in this research is called the Memetic Algorithm with Explicit Diversity Management (MA-EDM).
MA-EDM combines an effective local search and well-tested permutation-based genetic operators with a novel diversity-aware replacement strategy  that  aims  to  gradually  shift from  exploration  to exploitation  as  the  evolution progresses.
As shown in the experimental validation, quite a large number of new best-known solutions were attained for the most popular instances used in the literature.
This is quite an important achievement, especially considering that these benchmarks have been available for decades and that many different kinds of methodologies have been used to tackle them.
It is also worth noting that our proposal does not incorporate some of the problem-dependent knowledge about the structure of the problem that is taken into account in the most recent LOP solvers.
Thus, this research also shows that by applying proper general metaheuristic design principles, simpler and more effective LOP optimizers can be developed.

In the following section, the state of the art in the LOP area is summarized by introducing a chronological description of the most promising approaches and their main features.
Afterwards, our novel diversity-aware memetic algorithm is described in detail.
Finally, the experimental analysis is presented and the main conclusions are drawn. 
The experimental validation takes into account the best-performing metaheuristics proposed in the literature with the aim of showing that the advances are not only due to the fact that long-term executions are carried out.
The advantages of MA-EDM are quite remarkable.
Experiments to shed some light on the reasons behind the superiority of MA-EDM are also included.

\section{Literature review}\label{sec2}

A wide variety of approaches to deal with the LOP have been proposed~\cite{Marti:2011}, including exact methods like branch-and-bound or cutting planes, specific heuristics and multiple metaheuristic strategies.
In recent years, metaheuristics have gained more attention, especially for dealing with the most recent benchmarks, which include quite complex and large instances~\cite{Schiavinotto:2004, Ceberio:2015}.
The most relevant metaheuristics developed for the LOP are described briefly below.

In 1999, Laguna et al.~\cite{Laguna:1999} proposed a Tabu Search (TS) approach which incorporates strategies for search intensification and diversification. 
Two different definitions of basic neighborhoods and path relinking are introduced for search intensification, while the diversification relies on ad-hoc short-term and long-term components. 
The short-term component is focused on modifying items that have been changed infrequently. 
Meanwhile, the long-term component is used to perform restarts by generating solutions that are distant from the best solutions found so far.
This approach improved --- at the time --- the state-of-the-art heuristics~\cite{Chanas:1996}

As an alternative, Campos et al.~\cite{Campos:2001} proposed a Scatter
Search scheme that takes into account diversity in two of its components.
First, in order to generate initial diverse solutions, GRASP constructions were designed with a greedy function that combines attractiveness and diversity.
Second, the reference set selects half of its members by quality and half of its members by maximizing a measure of diversity. 
Regarding the improvement method, it is also based on the two neighborhoods proposed in~\cite{Laguna:1999}, while an ad-hoc
combination method based on votes is used.

Later, Schiavinotto and Stützle~\cite{Schiavinotto:2004} performed a search-space analysis of the available LOP instances. 
Depending on certain characteristics of the matrices, such as variability, number of zeros or symmetry, they realized that those based on real data are very different from random ones.
Moreover, the LOLIB instances~\cite{Grotschel:1984} were too small to pose a real challenge to state-of-the-art metaheuristic approaches, and even to exact algorithms, so an additional set of large, random, real-life like instances ---denoted as eXtended LOLIB (xLOLIB)--- was generated. 
The most complex instances present a low correlation between fitness and distance to optimal solutions, which suggests that 
exploring distant regions is important to performance.
As a result, an MA that includes an Iterated Local Search (ILS) for intensification and an adaptive restart strategy for diversification was proposed. 
For the local search, the neighborhoods were defined through the interchange or insert moves, and the swap operation of the current permutation. 
In addition to these neighborhoods, 
the sorting-reversing procedure introduced by Chanas and Kobylanski~\cite{Chanas:1996} was also applied.
Since the performance of an MA may depend strongly on the crossover operator, four different operators were tested. 
The most promising results were obtained with Cycle Crossover (CX) and Order Crossover (OB). 
The experimental study showed that the MA resulted in a much more robust performance; in particular, the restarting mechanism was key to the success.

Trajectory-based heuristics were further advanced by developing variants of Variable Neighborhood Search (VNS)~\cite{Garcia:2006}.
Several variants that include different intensification and diversification approaches were proposed.
Experimental analysis shows that the performance of VNS improves when combined with short-term memory TS and diversification schemes.
While proper solutions were attained in short-term executions, the solutions were not as promising as those output by MAs. 

%
In 2012, Martí, Reinelt and Duarte~\cite{Marti:2012} published an essential survey in the
field of the LOP. 
Their work offers a comprehensive study of the main existing heuristic and metaheuristic approaches for the LOP at the time.
For the heuristic methods, the authors studied some classical construction algorithms and
some improvement algorithms based on inserts ---moves or exchanges--- or on
local search.
The experimental study also included some multi-start methods~\cite{Chanas:1996}, as well as all the metaheuristics already introduced here: 
Tabu Search~\cite{Laguna:1999}, 
GRASP~\cite{Campos:2001}, 
Scatter Search~\cite{Campos:2001}, 
Genetic Algorithms~\cite{Schiavinotto:2004}, 
MA~\cite{Schiavinotto:2004}, 
ILS~\cite{Schiavinotto:2004}, and 
VNS~\cite{Garcia:2006}.
The experimental study also included a Simulated Annealing based on the noising algorithm by Charon and Hudry~\cite{Charon:2007, Charon:2010}.
For their analysis, the authors also compiled all the existing ---and widely used in related literature--- LOP benchmarks to release a new version of the LOLIB library. 
Their statistical analysis demonstrated that MA and ILS~\cite{Schiavinotto:2004} constitute the most promising approaches.

Subsequently, the potential of ILS was explored further.
Two general frameworks, based on the ILS and the Great Deluge Algorithm respectively, were proposed.
The main novelty relies on exploring the neighborhood in a more efficient way~\cite{Sakuraba:2010,Sakuraba:2015} with the aim of increasing the amount of iterations that might be evolved.
%

An alternative to reduce the computational cost of local search-based approaches was proposed by 
Ceberio, Men\-di\-buru, and Lozano~\cite{Ceberio:2015}. 
They conducted a theoretical analysis of neighborhoods that allowed them to identify specific positions where the indexes might generate local optima solutions.
This way, a restricted insert neighborhood is considered in conjunction with the best performing algorithms~\cite{Marti:2012}, MA and ILS~\cite{Schiavinotto:2004}, to generate the MAr and ILSr methods. 
The comparisons of MAr and ILSr ---at the level of solution quality but especially in terms of execution times---  demonstrated the superiority of these approaches.
In order to go further in their testing, the authors generated an extra set of more complex instances, known as xLOLIB2,
that were also integrated in LOLIB.
xLOLIB2 is currently the set of largest instances that is most discussed in the related literature.

In order to further improve the exploration behavior of ILS and MA, the Hybrid Exploration
Algorithm (HEA) was devised, which is a population-based metaheuristics based on ILS~\cite{Garcia:2019}.
Similarly to the proposal put forth in this paper, HEA moves the balance from exploration to exploitation during the execution.
This is done through a perturbation with a dynamic strength.
Note that differently to our proposal, this adaptation does not take into account the stopping criterion or elapsed 
period.
Additionally, a post-optimizer based on exactly solving submatrices of the instance by dealing with the LOP as a binary linear programming problem is applied. 
New best-known solutions were set to 77 out of the 78 instances from the xLOLIB~\cite{Schiavinotto:2004}.
%

Finally, with the aim of exploiting some structural properties of the LOP, Santucci and
Ceberio~\cite{Santucci:2020} designed a new approach, denoted by CD-RVNS, which involves the following components:
a VNS algorithm based on restricted neighborhoods, 
a randomized heuristic construction procedure that builds-up an LOP solution precedence by precedence, and 
a destruction procedure that removes precedences from an LOP solution with the aim of producing a new starting point for the construction procedure.
VNS considers two different neighborhoods, \textit{insert} and \textit{interchange}, by switching between them when the search is trapped in a local optimum. 
During the construction and destruction procedures, an LOP solution is represented as a set of pairwise precedences, thus differing from the classical linear representation of permutations, which is indeed used in the VNS phase. 
Dealing with both representations requires the implementation of conversion methods between representations. 
Although representing a linear ordering as a set of precedences requires more memory than encoding a linear permutation, such a representation allows constructive heuristics to smoothly build up the solutions, precedence by precedence, thus avoiding drastic changes in the objective value.

In summary, it is important to remark that most proposals call attention to the importance of diversification with the aim of avoiding stagnation.
However, recent design principles based on relating the diversity to the stopping criterion and elapsed period have never been applied in LOP solvers. 
Additionally, most recent strategies have focused on employing specific features of the structure of the LOP and/or hybridizing metaheuristics with exact solvers to advance the results further, meaning that the complexity of solvers 
has increased drastically, with algorithms containing several thousand lines of code.
Finally, trajectory-based approaches have generally been more successful in short-term executions, whereas in the long term, population-based approaches are usually more convenient.
While MAr reports most of the current best-known solutions, the trajectory-based ILSr and CD-RVNS solvers, currently
offer some of the best-known solutions for several instances, so when validating new solvers, it is important to 
consider both kinds of metaheuristics.

\section{Memetic Algorithm with Explicit Diversity Management: MA-EDM}\label{sec3}

Knowing near-optimal solutions for large instances is important for the analysis of the performance of
LOP solvers.
Given the hypothesis previously stated, the research is focused on developing a LOP solver that,
in the long-term, is able to find new best-known values for LOP instances of high complexity.
Diversity-aware strategies that gradually shift from exploration to exploitation 
as the evolution progresses have excelled in combinatorial problems~\cite{Segura2016}.
Thus, in order to attain our aim, a novel diversity-aware LOP solver, the Memetic Algorithm with Explicit Diversity Management (MA-EDM) is proposed.
MA-EDM was designed by considering some genetic operators that have already been proven effective 
for the LOP.
Algorithm \ref{alg:MADM} shows the pseudocode of our proposal.
It is a first-generation MA, so a static improvement individual strategy is applied.
Each component is described below in detail.

\begin{algorithm}[t]
\caption{Memetic Algorithm with Explicit Diversity Management}
\label{alg:MADM}
\begin{algorithmic}[1]
\Require {\textit{LOP instance}, \textit{N}(size of population), crossover operator (CX or OB), \textit{Stopping criterion (time)}}
\State{\textbf{Initialization:} Generate an initial population $P_0$ with $N$ individuals. Assign $i=0$.}
\State{\textbf{Local Search:} Perform a local search for every individual in $P_0$.} 
\State{\textbf{Diversity Initialization:} Calculate the initial desired minimum distance ($D_0$) as the mean distance among individuals in $P_0$.}
\While {the execution time has not finished}
	\State{\textbf{Mating Selection:} Perform binary tournament selection on $P_i$ in order to fill the mating pool with $N$ parents.}
	\State{\textbf{Variation:} Apply the crossover in the mating pool to create the set $O_i$ with $N$ offspring.}
	\State{\textbf{Local Search:} Perform a local search for each individual in $O_i$.}
	\State{\textbf{Survivor Selection:} Apply the replacement technique (BNP) to create $P_{i+1}$ by considering $P_i$ and $O_i$ as input.}
	\State{$i=i+1$}
\EndWhile
\State{\textbf{Return} best evaluated permutation.}
\end{algorithmic}
\end{algorithm}

Our approach starts by initializing a population with $N$ individuals randomly (line~1).
As stated previously, in most effective LOP solvers, solutions directly encode the permutation,
so this was the representation we selected.
The initialization considers a uniform random generator, meaning that each permutation is equiprobable.
Then, each solution is improved by applying a first-improvement stochastic hill-climber 
that considers the insert neighborhood (line~2)~\cite{Laguna:1999}. 
Note that this is one of the most successful neighborhood definitions, and neighbors are 
generated by simply moving a number of the given permutation to an alternative position while
shifting all the numbers in the intermediate positions.
Naive implementations of this local search evaluate the whole neighborhood in $O(N^3)$.
However, since the best insertion position of a number can be calculated in $O(N)$, neighbors are visited more systematically to reduce the total cost to
$O(N^2)$~\cite{Sakuraba:2010}.
Note that similar improvement strategies are applied in most state-of-the-art MAs~\cite{Sakuraba:2015,Santucci:2020}.
It is also important to note that the stochastic behavior of the hill-climber is maintained
because at each iteration, the order in which the numbers to move are considered is random.
Once a number is selected, if its best insertion move improves the current solution, 
it is accepted.
This is repeated until no neighbor improves the current solution.
Note that the local search with the complexity reduction technique allows achieving
local optima quite efficiently, even for the largest instances considered in our
experimental validation.
For example, in instances with $n = 1000$, the time to reach a local optimum, when starting from a randomly created solution, is less than one second.
This means that it can be directly integrated into MAs without the need to consider alternative stopping criteria.
 
An important feature of our approach is that it considers diversity explicitly. 
This is done with the replacement strategy, which works by setting a minimum
desired distance that is updated during the run.
The details are given later, but similarly to~\cite{Hernandez:21}, the distances that appear in the
initial population are used to set the initial desired distance ($D_0$) (line~3). 
Note that this has the advantage of adapting the execution automatically to the 
given instance~\cite{Hernandez:21}, in the sense that usually, for the most complex instances where
a higher diversity is required to reach high-quality solutions, the improved initial population
tends to exhibit larger distances.
$D_0$ is calculated as the mean distance among all the individuals in the initial population.
In order to calculate $D_0$, a metric or a distance-like function must be selected.
In this case, the permutation deviation distance ($D_{dev}$) \cite{sevaux2005permutation} is employed. 
The dissimilarity between two permutations $x$ and $y$ with $p$ elements is defined as follows:

$$ D_{dev}(x, y) =  \sum_{i=1}^{p}{ \| i - pos( x[i], y ) \| }  $$

where $pos(v, perm[])$ is the index of the value $v$ in the permutation $perm$.

\vspace{0.25cm}
The \textit{permutation deviation distance} was selected because the position of a number 
establishes the amount of values that are summed from the associated column.
Thus, summing the absolute distance between the positions for each number is meaningful.
Note that this metric has a broader range than other typical metrics applied to permutation-based representations, such as the exact match distance, which provided important advantages for the Job Shop Scheduling problem~\cite{Hernandez:21}.
Moreover, it can be calculated in $O(N)$, meaning that the impact on 
performance is low, which also justifies
its selection.
Any of the metrics described in~\cite{sevaux2005permutation} could also be used, and some of them might even yield better results.
However, given the advantages of the permutation deviation distance mentioned previously, and that this work considers a large amount of complex instances that require computationally expensive experiments, the analyses of alternative metrics is left for future work.

As in standard MAs, MA-EDM evolves a set of generations until a given stopping criterion is reached (line~4-10).
In order to provide fair comparisons, and since very different 
types of proposals, including trajectory-based and population-based approaches,
are considered in our validation, the stopping criterion was set by time.
At each generation, a set of $N$ parents is selected using binary tournaments (line~5).
Then, the variation phase is applied (line~6).
In this case, the variation consists only of a crossover operator, and the mutation
is omitted.
The logic behind this decision is that in many cases, MAs achieve proper solutions
even if mutation is not used~\cite{Neri:11}.
The local search performs modifications that allow new alleles to appear, 
so the increased complexity of including an additional operator is not always justified~\cite{Tsai:04}.
This is especially the case when applying diversity-aware strategies~\cite{Hernandez:21}, because one of the other roles of mutation is to avoid the premature loss of diversity, which in our case
is preserved explicitly in the replacement phase.
Regarding the crossover operators, two of the most effective are the 
Cycle Crossover (CX) \cite{Merz2000} and the Order Crossover (OB) \cite{Davis1991}.
A clear winner is not identified in the related literature~\cite{Schiavinotto:2004}, so both of them are considered in this research.
The CX operator is based on the principle of inheriting values in its exact positions.
Thus, each number in the offspring permutations shares its position with one of the two parents.
First, both offspring inherit the elements that share the same position in both
parents, by just copying them in their original position.
Then, the remaining elements are used to establish cycles that fulfill the following property: consecutive elements in the cycles share the same position in the parents.
Note that with this property, the cycles can be identified unequivocally.
For each number, one of the parents is selected randomly and the elements belonging to the associated cycle are copied in the same positions as in the selected parent. 
The OB operator is concerned with preserving the partial order among elements.
First, one of the parents is cloned.
Then, half of the positions are selected randomly, and they are sorted according to the order appearing in the second parent~\cite{Davis1991}.
Note that in order to generate the second offspring, the roles of the parents are interchanged.
In order to improve the generated offspring, they are subjected to the same improvement strategy that was applied to the initial population (line~7).

Finally, the survivor selection strategy is applied to select the members of the next population using the Best Non-Penalized (BNP) survivor selection strategy (line~8).
The main design principle behind BNP is similar to the
one discussed in~\cite{Segura2016} and relies on introducing a dependency between the selection of survivors and the
ratio between the elapsed execution period and the stopping criterion.
Specifically, the aim is to promote a larger degree of exploration in the
initial phases of the execution, and then to gradually shift towards exploitation as the
end of the optimization process nears.
Note that BNP was successfully applied in the optimization of the Job Shop Scheduling Problem~\cite{Hernandez:21}.
However, to our knowledge, it has never been applied in the Linear Ordering Problem.

\begin{algorithm}[!ht]
\caption{Best Non-Penalized (BNP) Survivor Selection Strategy}
\label{alg:BNP}
\begin{algorithmic}[1]
\Require {Population $P=\{p^1,...,p^N\}$, Offspring $ O = \{o^{1}, ..., o^{N}\}$}
\State{$D=D_0 - \frac{T_{Current}}{T_{Total}} \times D_0$}
\State{$S_{Eligible} = P \cup O$} 
\State{$S_{Penalized}=\emptyset$}
\State{$S_{NewPopulation}=\emptyset$}
\While {$\vert S_{NewPopulation}\vert < N$}
\For{$individual \in S_{Eligible}$}
    \State{$individual.dci = DCI(individual, S_{NewPopulation})$}
        \If{$individual.dci < D$}
            \State{$S_{Eligible} = S_{Eligible} \setminus{\{individual\}}$}
            \State{$S_{Penalized} = S_{Penalized} \cup \{individual\} $}
        \EndIf
\EndFor
    \If{$S_{Eligible} \neq \emptyset$}
        
        \State {$new\_selected = \argmax_{x\in S_{Eligible}} {x.obj}$}// Ties are broken randomly
    \Else
         \For{$individual \in S_{Penalized}$}
            \State{$individual.dci = DCI(individual, S_{NewPopulation})$}
         \EndFor
         \State {$new\_selected = \argmax_{x\in S_{Penalized}} {x.dci}$}// Ties are broken randomly
    \EndIf
    \State{$S_{NewPopulation} = S_{NewPopulation} \cup \{new\_selected\}$}
\EndWhile
\State \Return $S_{NewPopulation}$
\end{algorithmic}
\end{algorithm}

Algorithm~\ref{alg:BNP} shows the pseudocode of BNP.
It receives as its input the current population and offspring, and it outputs the population of the next generation.
BNP is based on promoting both quality and diversity, and the second aim is attained
by avoiding the selection of individuals that are too close to each other.
The definition of too close is dynamic.
Specifically, a threshold distance ($D$) that is used to identify too-close individuals is calculated with the formula given in line~1.
In this formula, $T_{current}$ is the number of elapsed seconds from the start of
the optimization process, and $T_{total}$ is the stopping criterion (also in seconds).
Thus, the initial distance ($D_0$) previously discussed is reduced linearly, so that at the end of the optimization process, $D$ is set to zero.
This means that, as the evolution progresses, the promotion of diversity is reduced and, consequently, more importance is given to the quality.
This transformation is performed gradually, resulting in a steady shift from exploration to intensification.
BNP manages three sets of individuals.
The eligible set ($S_{Eligible}$) contains individuals that have not yet been selected to survive and that are not too close (distance lower than $D$) to any of the solutions already selected to survive.
The penalized set ($S_{Penalized}$) contains individuals that have not yet been selected to survive and whose distance to a survivor is lower than $D$.
Finally, the new population set ($S_{NewPopulation})$  contains the survivors.
Initially, $S_{Eligible}$ is filled with the union of the current population and offspring (line~2), whereas $S_{Penalized}$ and $S_{NewPopulation}$ are empty (lines~3 and 4).
Finally, the process iteratively selects survivors until $S_{NewPopulation}$ is filled with $N$ individuals (line~5 to 22).
At each iteration, the distance condition is checked for each individual in $S_{Eligible}$ (lines~6 and 7), and those that are too close to an already selected survivor are moved to $S_{Penalized}$ (line~8 to 11).
Note that in the pseudocode, $DCI$ stands for \textit{Distance to Closest Individual}, and is defined as follows: $DCI(x, P) = min\{ distance(x, y)\ \vert \ y \in P \}$, where $distance(x,y)$ refers to the deviation permutation distance already discussed.
In those iterations where $S_{Eligible}$ is not empty, i.e. there are sufficiently distant individuals, the best among them (the one that maximizes the objective function) is selected to survive (lines~13-14 and 21).
The opposite case represents situations where the diversity is not as large as desired.
Thus, the penalized solution with the longest distance to the closest survivor is selected (lines 16-21).
\begin{figure}[t]
\centering 
\subfloat[Initial phase]{\label{fig:bnp1}\includegraphics[width=0.30\textwidth]{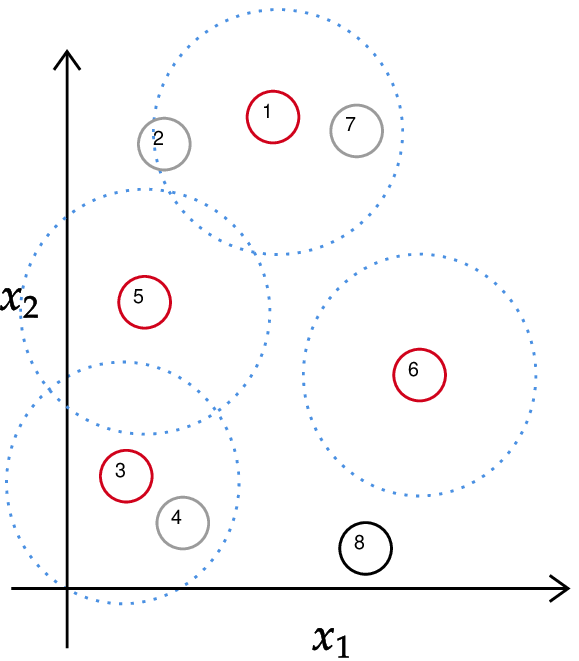}}
\vspace{0.1cm}
\subfloat[Intermediate phase]{\label{fig:bnp2}\includegraphics[width=0.30\textwidth]{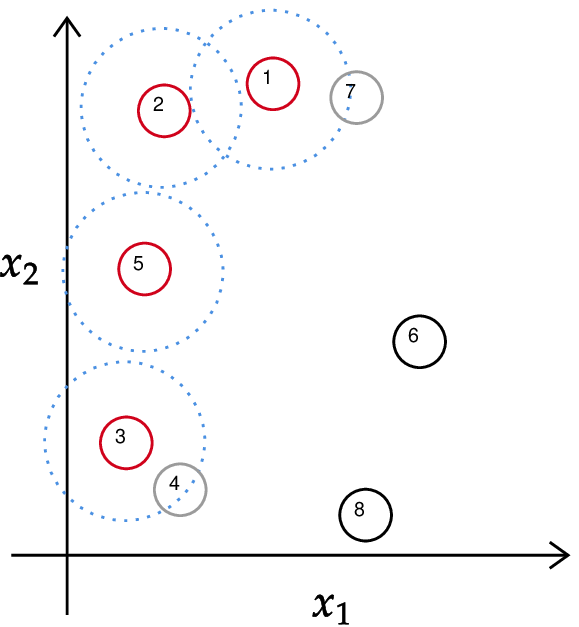}}
\vspace{0.1cm}
\subfloat[Final phase]{\label{fig:bnp3}\includegraphics[width=0.30\textwidth]{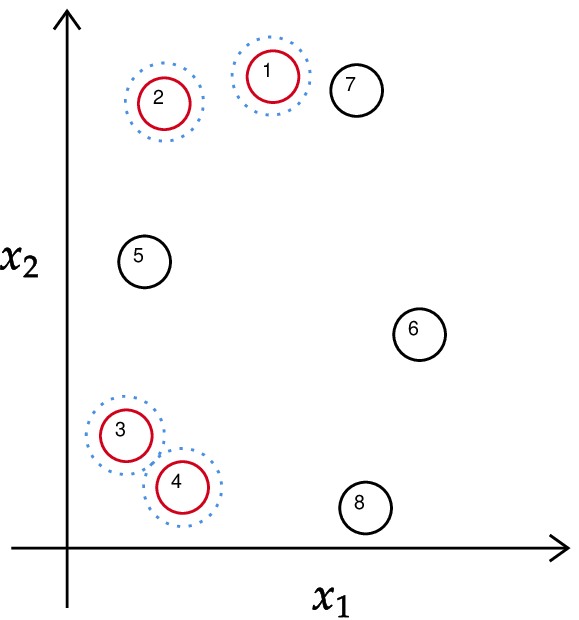}}
\qquad
\caption{Illustration of the BNP replacement process for a continuous optimization problem with two variables} 
\label{fig:bnp_process} 
\end{figure}

In order to further illustrate the replacement operator, Figure~\ref{fig:bnp_process} shows the behavior of BNP in three different scenarios: initial, intermediate and final phases of the run.
Particularly, in each case, it considers that its input is a set of eight solutions (union of previous population and offspring) from which four solutions must be picked up as survivors.
Each subfigure represents the status at the end of the BNP process.
Note that the threshold distance to perform the penalties is reduced over time, meaning that each of the previous scenarios consider different threshold distances.
In order to facilitate the illustration, this figure considers that BNP is applied to a continuous problem with two variables ($x_1$ and $x_2$).
Each small circle with a number represents a solution and the value inside the circle represents its relative quality, meaning the solution tagged with a $1$ is the best one, whereas the one tagged with a $8$ is the worst-one.
Solutions with a red circle are the survivors ($S_{NewPopulation}$), the ones with a gray circle are the penalized ones ($S_{Penalized}$), whereas the ones with a black circle are not penalized nor selected ($S_{Eligible}$).
Also note that the blue circle surrounding each survivor represents the region of the search space that is penalized because of such a selection.
Note that, at each iteration of BNP, the best non-penalized solution is selected until four or them are marked in red color.
Due to the different threshold sizes, the survivors marked in subfigure (a) are more focused on diversity than on quality, whereas subfigure (c) represents a case with a higher selection pressure.
In this particular example, it selects the best solutions even if it implies an important reduction of diversity.
One of the keys of BNP is precisely to perform a smooth reduction of the threshold distance, with the aim of shifting gradually from exploration to exploitation.

\section{Experimental Study}\label{sec4}
%
This section is devoted to validate the advances attained by MA-EDM. 
One of the main aims of this research is to obtain new best-known solutions for the most complex sets of instances, so that they can be used as new reference solutions.
However, it is also important to show that such solutions are not attainable by state-of-the-art solvers when similar computational resources are granted. Thus, the set of experiments aims to provide new best-known solutions, compare MA-EDM against state-of-the-art solvers and shed some light on the reasons behind the superior performance of MA-EDM.

In the LOP field, several sets of instances with different dimensions and degrees of complexity have been proposed and integrated in LOLIB. 
The experimental validation considers nine of the most popular benchmarks sets~\cite{Marti:2012}, which are the following: IO~\cite{Grotschel:1984}, MB~\cite{Mitchell:2000}, SGB~\cite{Knuth:1993}, RandB~\cite{Marti:2012}, RandA1~\cite{Laguna:1999}, RandA2~\cite{Campos:2001}, Spec~\cite{Marti:2012}, xLOLIB~\cite{Schiavinotto:2004} and xLOLIB2~\cite{Ceberio:2015}.
Some of them, especially the four first sets of benchmarks, are relatively easy and are almost systematically solved by current solvers.
Thus, they are useful just to check the robustness of new proposals.
The rest of the benchmark sets contain some larger and more difficult cases, and even though some of them were proposed more than 30 years ago, they have yet not been solved to optimality.
Thus, they are normally used to check new advances in the field.
\begin{table}[]
\centering
\begin{tabular}{cccc}
\textbf{Set} & \textbf{Instances} & \textbf{Size(\#Instances)}                                                      & \textbf{Optimum Known} \\ \hline
IO           & 50                 & \begin{tabular}[c]{@{}c@{}}44(31)\\ 50(4)\\ 56(11)\\ 60(3)\\ 79(1)\end{tabular} & 50              \\ \hline
MB           & 30                 & \begin{tabular}[c]{@{}c@{}}100(5)\\ 150(10)\\ 200(10)\\ 250(5)\end{tabular}     & 30              \\ \hline
SGB          & 25                 & 75(25)                                                                          & 25              \\ \hline
RandB        & 90                 & \begin{tabular}[c]{@{}c@{}}40(20)\\ 44(50)\\ 50(20)\end{tabular}                & 71              \\ \hline
RandA1       & 100                & \begin{tabular}[c]{@{}c@{}}100(25)\\ 150(25)\\ 200(25)\\ 500(25)\end{tabular}   & -             \\ \hline
RandA2       & 75                 & \begin{tabular}[c]{@{}c@{}}100(25)\\ 150(25)\\ 200(25)\end{tabular}             & 25              \\ \hline
Spec         & 37                 & \begin{tabular}[c]{@{}c@{}}11-36(8)\\ 50-452(29)\end{tabular}                   & 30              \\ \hline
xLOLIB       & 78                 & \begin{tabular}[c]{@{}c@{}}150(39)\\ 250(39)\end{tabular}                       & -               \\ \hline
xLOLIB2      & 200                & \begin{tabular}[c]{@{}c@{}}300(50)\\ 500(50)\\ 750(50)\\ 1000(50)\end{tabular}  & -               \\ \hline
\end{tabular}
\caption{Summary information of the sets of instances that comprise LOLIB.}
\label{tab:summary-lolib}
\end{table}
Table~\ref{tab:summary-lolib} summarizes some of the most relevant features of the sets of instances included in LOLIB.
Particularly, it shows the number of instances that comprises each set, the sizes of the instances, and the number of instances where the current best-known solution matches the current best-known upper bound.
Note that only in the smallest instances, there is a match between the best-known solutions and the best-known upper bounds.
As it is shown below, the best-known solutions are improved in this paper for 293 instances.
However, the best-known upper bounds were not attained in any of those instances.
In fact, the distances between the upper bounds and the best-known solutions is quite large, so probably there is also room for improvement in terms of the upper bounds.

Regarding state-of-the-art solvers, the validation was conducted by considering both trajectory- and population-based metaheuristics.
For the trajectory-based methods, the ILSr~\cite{Ceberio:2015}
 and CD-RVNS~\cite{Santucci:2020} methods were used.
Note that CD-RVNS yields better results than ILSr in most cases.
However, there are some specific instances where ILSr provides a highly effective search, so it was maintained in our validation.
Regarding the population-based alternatives, MAr~\cite{Schiavinotto:2004} is the current clear winner, and similarly to MA-EDM, it is also a Memetic Algorithm, so it was selected to perform the validation.
Note that the source codes of all these methods are publicly available.
In order to facilitate the replicability of our experiments, our source code is also freely available~\footnote{\url{https://github.com/carlossegurag/LOP_MA-EDM}. In addition to the source code, the datasets generated and/or analyzed during the current study are available in this repository in the RawData directory. Moreover, this repository contains the best-known solutions generated for each instance in the BKS.zip file.}.

Three kinds of experiments are presented below.
First, the tuning step of MA-EDM is described.
Second, MA-EDM is compared against the selected state-of-the-art solvers.
Finally, some analyses that provide a better understanding of the differences between the
dynamics of the populations of MAr and MA-EDM are detailed.
Note also that in order to provide sound conclusions, and taking into account the stochastic behavior of the solvers involved, all our executions were repeated 30 times and proper statistical tests were applied.
Specifically, the results were compared with the \textit{scmamp library}~\cite{Calvo2016} and Scikit-learn~\cite{Pedregosa:11}, following some of the guidelines described in~\cite{Garcia:10}.
All the tests were applied assuming a significance level of 5\%.

\subsection{Parameterization}\label{subsec1}

The parameterization of MA-EDM is a relatively simple task.
The two only decisions that had to be made to perform the experimentation were the selection of the crossover operator and the population size.
As mentioned earlier, regarding the crossover operator, it was clear that CX and OB had yielded the best overall results~\cite{Schiavinotto:2004},
so both of them were considered in the tuning phase.
Regarding the population size, it is well known that adequate values usually depend on other components of the optimizer, such as the stopping criterion, in non-trivial ways~\cite{Eiben:15,Zhang:10}.
Since diversity-aware techniques usually require long executions and they alter the dynamics of the population, using previously reported experiments to select the population size is not too meaningful, so we decided to test the following values: 25, 100, 200, 300, 400, 500 and 600.

In order to understand the implications of the parameterization on the results, four instances of the benchmark with the largest instances (xLOLIB2) were considered.
Since xLOLIB2 consists of instances with four different sizes, one instance of each size was selected.
They were \textit{Nbe75eec300}, \textit{Nbe75eec500}, \textit{Nbe75eec750}, and \textit{Nbe75eec1000}.
The fourteen different configurations resulting from combining two crossover operators and seven population sizes were executed 30 times for each instance.
The stopping criterion was set to 4 hours.

\begin{figure}[t]
\centering 
\subfloat[Nbe75eec300]{\label{fig:1}\includegraphics[width=0.48\textwidth]{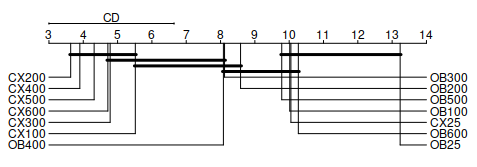}}
\subfloat[Nbe75eec500]{\label{fig:2}\includegraphics[width=0.48\textwidth]{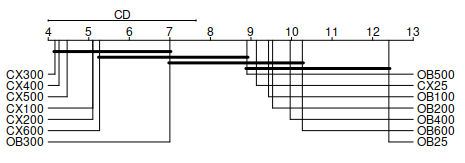}}
\qquad
\subfloat[Nbe75eec750]{\label{fig:3}\includegraphics[width=0.48\textwidth]{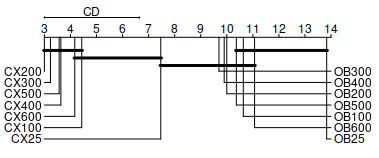}}
\subfloat[Nbe75eec1000]{\label{fig:4}\includegraphics[width=0.48\textwidth]{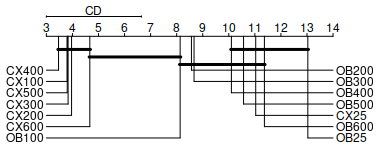}}
\qquad
\caption{Critical differences plot for each instance selected, considering the 14 configurations applied in the tuning phase} 
\label{fig:configuration} 
\end{figure}

The results of each instance were analyzed independently.
First, the Kruskal Wallis test was applied as a omnibus test.
In every case, significant differences appeared so Mann Whitney with Hommel's correction was used
to identify the pair-wise differences~\cite{del2013srcs}.
Additionally, with the aim of visualizing the results, the Nemenyi test was applied. 
The critical difference plot of each instance is shown in Figure~\ref{fig:configuration}.
Note that in this kind of plot, the average rank position of each model is shown and horizontal lines join schemes with no statistically significant differences.
Each scheme is denoted with a string formed by the concatenation of the crossover operator and the population size.
The advantage of the CX operator is quite clear.
In the four instances, at least the six best configurations use the CX operator.
Moreover, there was only one instance (\textit{Nbe75eec500}) where the results attained with the best configuration that relies on the OB operator was not statistically different than the ones reported by the best CX operator.
Thus, the superiority of the configurations with the CX operator is considerable.
As for the population sizes, the conclusions are not as clear.
First, an important asset of MA-EDM is that for the configurations with the CX operators, most of the configurations with different population sizes appear linked, meaning that results are not very sensitive to the population size.
Only the CX25 configuration is significantly different in some of the instances.
Thus, there is quite a large range of admissible values.
Since CX200 was the best ranked in two instances, and its ranking difference was never larger than half of the critical difference, this was the configuration selected to continue our validation.
Note that selecting any of the other configurations in the 100 to 600 range would probably lead to similar conclusions.

Finally, note that while the Nemenyi critical difference plot is quite illustrative, 
the more advanced statistical tests Kruskal Wallis (as omnibus) and Mann Whitney test with the Hommel's correction (as post-hoc) are more advisable~\cite{del2013srcs}.
When comparing CX200 to the remaining configurations with these tests, the same conclusions than the ones obtained with the Nemenyi test could be drawn.

\subsection{Validation against the state of the art}\label{subsec2}

One of the most important objectives of our research is to design a proper long-term solver, i.e.,
an optimizer that is able to profit from long-term executions by avoiding stagnation problems.
In order to validate this aim, the three state-of-the-art optimizers selected and the MA-EDM
were executed with every instance of the nine benchmark sets by setting the stopping criterion
to 4 hours.
Table~\ref{tab:parameterization} shows the parameterization applied in each optimizer.
Note that this is the parameterization recommended by the authors, but in the case of MAr, a parameterization
study similar to the one presented above was carried out to set the population size, and the value $N = 200$ was selected.
Each execution was repeated 30 times with different seeds.

\begin{table}[]
\centering
\begin{tabular}{cc}
\hline
\textbf{Optimizer} & \textbf{Parameterization}  \\ \hline
ILSr               & $\epsilon=0.0001$, $i\_{change}=7$, $n\_{ni}=750$   \\
MAr                & $N=200$                      \\
CD-RVNS                & Parameter-less Algorithm   \\
MA-EDM             & N=200, Crossover Operator=CX \\ \hline
\end{tabular}
\caption{Parameterization applied in each state-of-the-art solver and MA-EDM}
\label{tab:parameterization}
\end{table}

In order to summarize our results, statistical tests similar to those from the previous experiment were applied.
Specifically, the results attained for each instance with the four optimizers tested were compared by using the Kruskal Wallis test as an Omnibus test, followed by the Mann Whitney with the Hommel's correction to detect pair-wise differences.
Note that, for every instance, MA-EDM always attained an equal or larger mean and median than the
rest of the methods.
Table~\ref{tab:summary of excecutions} shows, for each benchmark set, the number of instances where 
the differences between MA-EDM and each state-of-the-art method selected were statistically significant
(column $>$) and not statistically significant (column =).
As we can note, the amount of instances where the results attained by MA-EDM are statistically superior is quite noticeable.
Particularly, the superiority appeared in 439, 289 and 420 instances out of the 685 instances, when comparing MA-EDM to ILSr, MAr and CD-RVNS, respectively.
In the smallest benchmarks, no much statistically significant differences appeared.
However, as the complexity grows, the number of instances with statistically significant differences in
favor of MA-EDM clearly increases.
Particularly noticeable is the case of xLOLIB2, where significant differences appeared in the whole
benchmark set.

\begin{table}[t]
\centering
	\resizebox{12cm}{!} {
\begin{tabular}{ccccccccccccccccccc}
\hline
\multirow{2}{*}{} & \multicolumn{2}{c}{IO} & \multicolumn{2}{c}{MB} & \multicolumn{2}{c}{RandA1} & \multicolumn{2}{c}{RandA2} & \multicolumn{2}{c}{RandB} & \multicolumn{2}{c}{SGB} & \multicolumn{2}{c}{Spec} & \multicolumn{2}{c}{xLOLIB} & \multicolumn{2}{c}{xLOLIB2} \\ \cline{2-19} 
 & = & \textgreater{} & = & \textgreater{} & = & \textgreater{} & = & \textgreater{} & = & \textgreater{} & = & \textgreater{} & = & \textgreater{} & = & \textgreater{} & = & \textgreater{} \\ \hline
MA-EDM vs ILSr & 50 & 0 & 29 & 1 & 2 & 98 & 42 & 33 & 90 & 0 & 0 & 25 & 33 & 4 & 0 & 78 & 0 & 200 \\
MA-EDM vs MAr & 50 & 0 & 30 & 0 & 67 & 33 & 75 & 0 & 90 & 0 & 25 & 0 & 36 & 1 & 23 & 55 & 0 & 200 \\
MA-EDM vs CD-RVNS & 50 & 0 & 30 & 0 & 13 & 87 & 51 & 24 & 90 & 0 & 0 & 25 & 31 & 6 & 0 & 78 & 0 & 200 \\ \hline
Total instances & \multicolumn{2}{c}{50} & \multicolumn{2}{c}{30} & \multicolumn{2}{c}{100} & \multicolumn{2}{c}{75} & \multicolumn{2}{c}{90} & \multicolumn{2}{c}{25} & \multicolumn{2}{c}{37} & \multicolumn{2}{c}{78} & \multicolumn{2}{c}{200} \\ \hline
\end{tabular}
}
\caption{Summary of statistical tests: column $>$ shows the number of instances where statistically significant differences appeared in favor of MA-EDM, and column $=$ counts the cases with no significant differences}
\label{tab:summary of excecutions}
\end{table}

\begin{table}[t]
\centering
	\resizebox{12cm}{!} {
\begin{tabular}{cccccc}
\hline
\multicolumn{1}{c}{\textbf{}} & \textbf{BKS} & \multicolumn{1}{c}{\textbf{ILSr}} & \textbf{MAr} & \multicolumn{1}{c}{\textbf{CD-RVNS}} & \textbf{MA-EDM} \\ \hline
Ntiw56n54\_300 & 2654641 & 2657882 & 2661516 & 2650252 & \textbf{26629001} \\
Nt75e11xx\_300 & 145044286 & 145093367 & 145297161 & 144958876 & \textbf{145388484} \\
Nt65f11xx\_500 & 29285209 & 29284570 & 29336499 & 29247869 & \textbf{29365313} \\
Nt75e11xx\_500 & 374756997 & 374780909 & 375603945 & 374381639 & \textbf{375998123} \\
Ntiw56n58\_750 & 20633749 & 20661491 & 20704523 & 20645181 & \textbf{20720453} \\
Nt65b11xx\_750 & 129443488 & 129523720 & 129868029 & 129349693 & \textbf{130011100} \\
Nbe75eec\_1000 & 122183020 & 122395155 & 122743873 & 122255253 & \textbf{122852580} \\
Nt70f11xx\_1000 & 177714930 & 177754220 & 178289425 & 177802295 & \textbf{178387875} \\ \hline
\end{tabular}
}
\caption{Mean objective function attained in eight selected instances from xLOLIB2}
\label{tab:summary of eight instances xLOLIB2}
\end{table}

Since xLOLIB2 contains the largest instances, eight instances were selected from this benchmark to show the results in more detail.
They were two instances from each of the different sizes present in the benchmark set.
Table~\ref{tab:summary of eight instances xLOLIB2} shows the mean objective function attained by each method, as well as the current \textit{Best-Known Solution} (BKS).
The improvement provided by MA-EDM is noticeable, and the differences with respect to the BKS are quite significant, which indicates that current solvers are not as effective with these large and complex instances, and it is precisely in the most complex cases where MA-EDM excels. 
Table \ref{tab:best_knowns_per_algorithm} shows the amount of new BKS found by each method in the different benchmark sets.
Note that when more than one method found the same new BKS, it is counted in all of them.
The advantages provided by MA-EDM are quite clear.
In fact, there was no instance where MA-EDM was exceeded in terms of the best or mean attained value.

In order to facilitate future validations, the supplementary material attached to this work contains the detailed results for the 685 instances considered over the course of this research.
For each instance, the best and mean solution attained by each method is shown.
In addition, the BKS is included.
MA-EDM was able to set new bounds in 293 cases out of the 685 instances, and it was always able to attain the BKS, meaning that the increase in robustness is huge in comparison to current solvers.

\begin{table}[]
\centering
\resizebox{\textwidth}{!}{%
\begin{tabular}{cccccc}
\hline
\multicolumn{6}{c}{\textbf{New BKS values per sets of instances}}                                                                  \\ \hline
\textbf{Instance/Algorithm} & \textbf{ILSr} & \textbf{MAr} & \textbf{CD-RVNS} & \textbf{MA-EDM} & \textbf{\# Instances Improved} \\ \hline
\textbf{IO}                 & 0             & 0            & 0            & 0               & 0                           \\
\textbf{MB}                 & 0             & 0            & 0            & 0               & 0                           \\
\textbf{SGB}                & 0             & 0            & 0            & 0               & 0                           \\
\textbf{Spec}               & 3             & 3            & 1            & 3               & 3                           \\
\textbf{RandB}              & 0             & 0            & 0            & 0               & 0                           \\
\textbf{RandA1}             & 2             & 10           & 4            & 38              & 38                          \\
\textbf{RandA2}             & 1             & 4            & 3            & 4               & 4                           \\
\textbf{xLOLIB}             & 0             & 16           & 5            & 48              & 48                          \\
\textbf{xLOLIB2}            & 0             & 0            & 0            & 200             & 200                         \\ \hline
\textbf{Total of new BKS}   & 6             & 33           & 13           & 293             & \textbf{293}                \\ \hline
\end{tabular}%
}
\caption{New BKS found in the experimentation process}
\label{tab:best_knowns_per_algorithm}
\end{table}

\subsection{Analysis of Diversity}\label{subsec3}

\begin{figure}[t] 
\centering 
\subfloat[Nbe75eec300]{\label{fig:9}\includegraphics[width=0.48\textwidth]{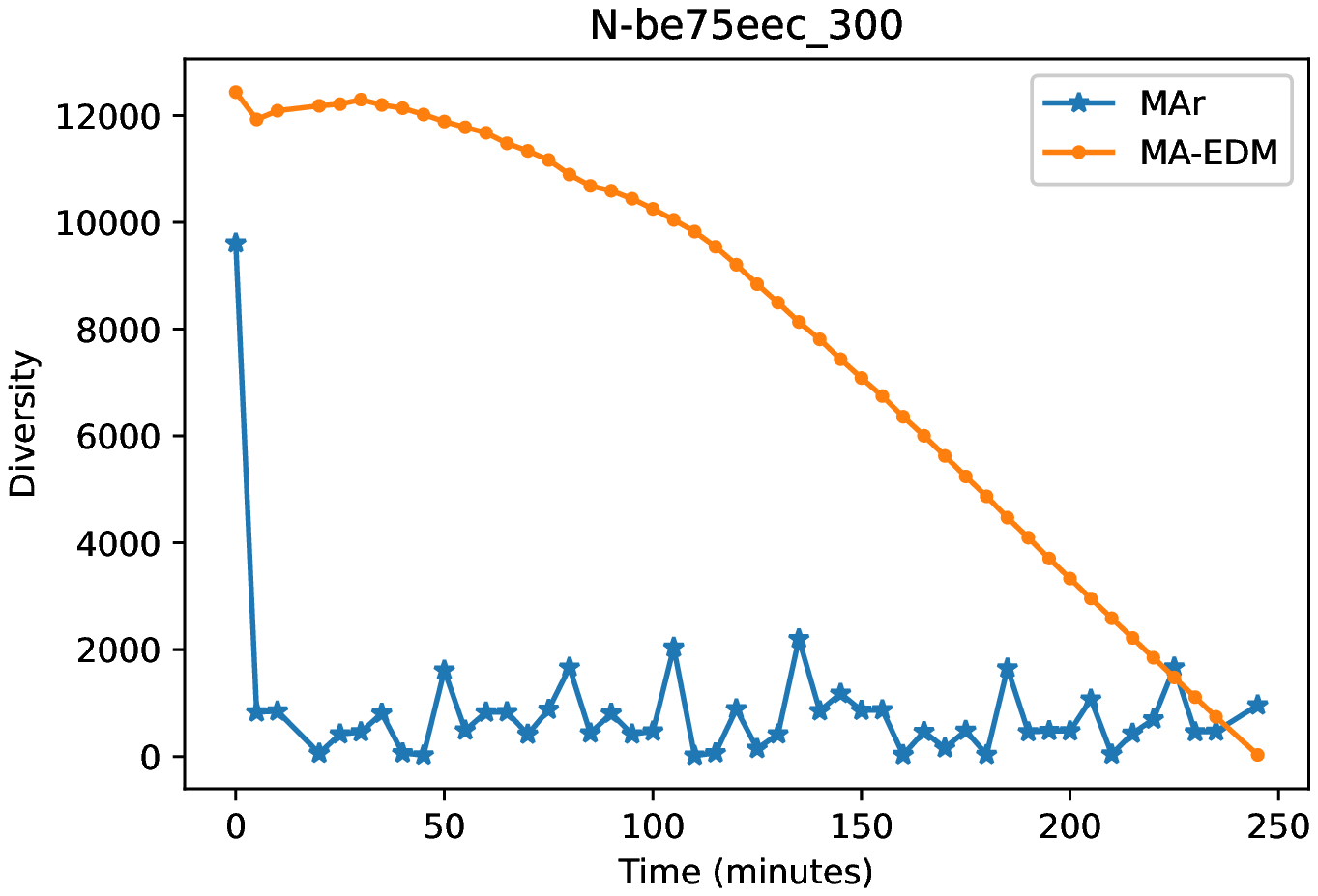}}
\subfloat[Nbe75eec500]{\label{fig:10}\includegraphics[width=0.48\textwidth]{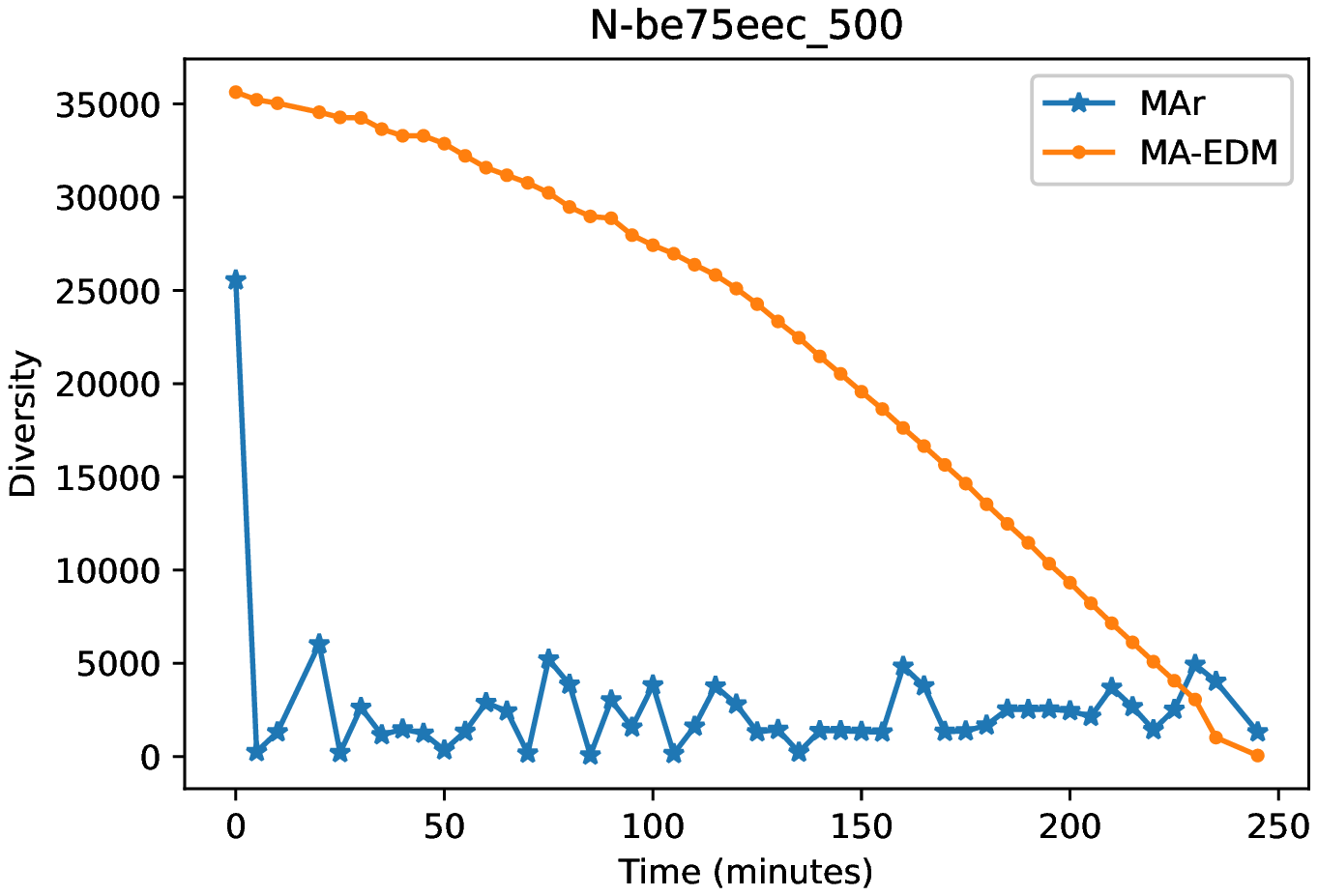}}
\qquad
\subfloat[Nbe75eec750]{\label{fig:11}\includegraphics[width=0.48\textwidth]{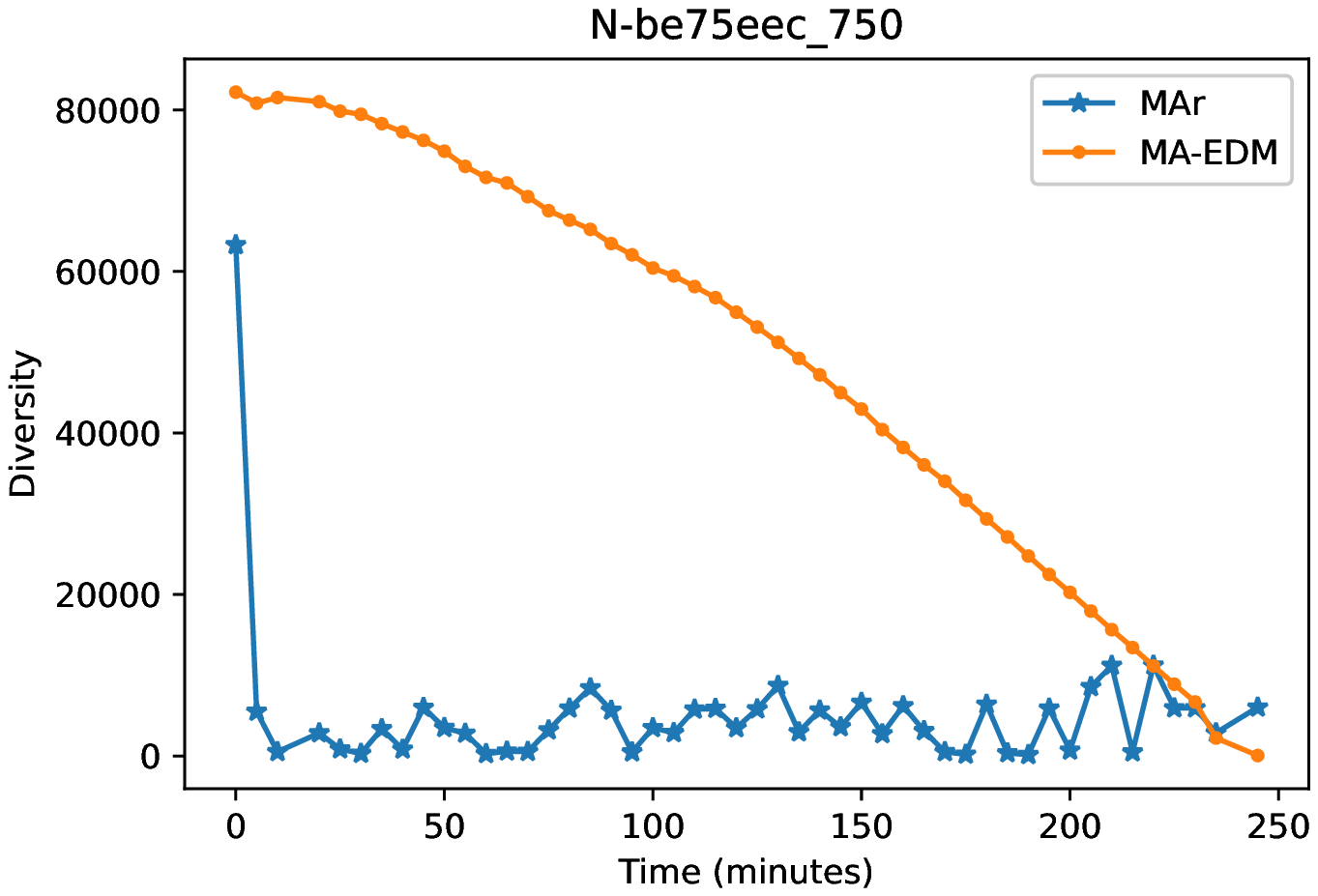}}
\subfloat[Nbe75eec1000]{\label{fig:12}\includegraphics[width=0.48\textwidth]{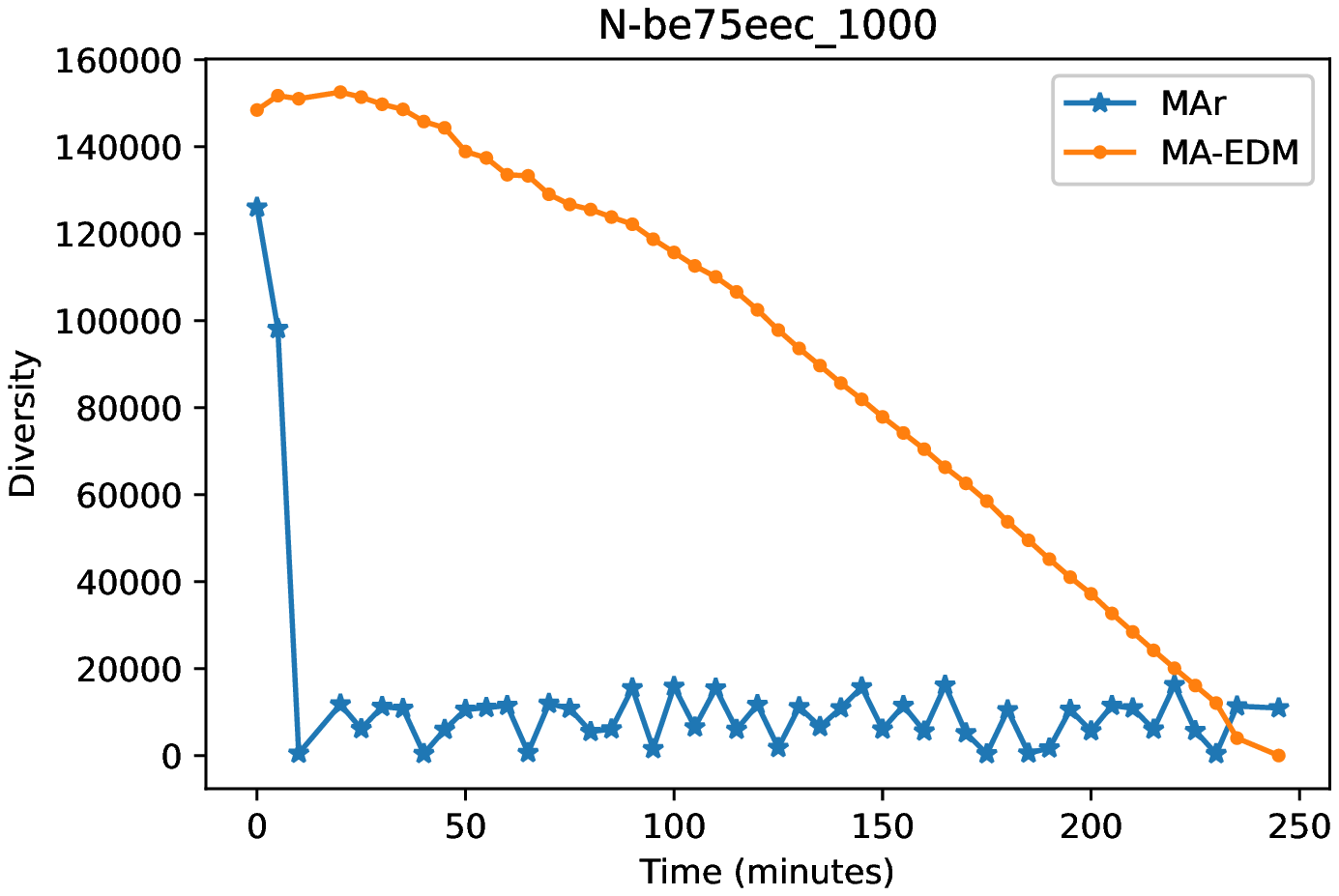}}
\qquad
\caption{Evolution of the diversity (mean permutation deviation distance) during the run in four selected instances (mean of 30 executions)} 
\label{fig:diversity_behavior} 
\end{figure}

One of  the most important distinguishing features of MA-EDM is the explicit
management of the diversity present in the population.
In order to understand the different dynamics of the populations of MA-EDM and MAr,
four instances were selected and the intermediate results were saved.
Figure~\ref{fig:diversity_behavior} shows the trend in the population's diversity
during the run.
Specifically, the diversity was calculated as the mean of the permutation deviation 
distances appearing in the population, and the results of 30 independent executions were averaged.
The four instances tested showed the same trend.
In MA-EDM, the reduction in diversity was almost linear, which
is explained by the linear reduction in the desired distance that is considered in the replacement phase.
By contrast, MAr reduced its diversity quite quickly, maintaining a low but non-zero diversity.
This means that MAr performs the search by placing all the individuals in a relatively narrow region and all the individuals are moved relatively close
to one another.
In comparison, MA-EDM initially maintains the individuals quite distant from one another, and as the execution progresses, they are gathered in the most promising regions located so far.

\begin{figure}[t] 
\centering 
\subfloat[Nbe75eec300]{\label{fig:5}\includegraphics[width=0.48\textwidth]{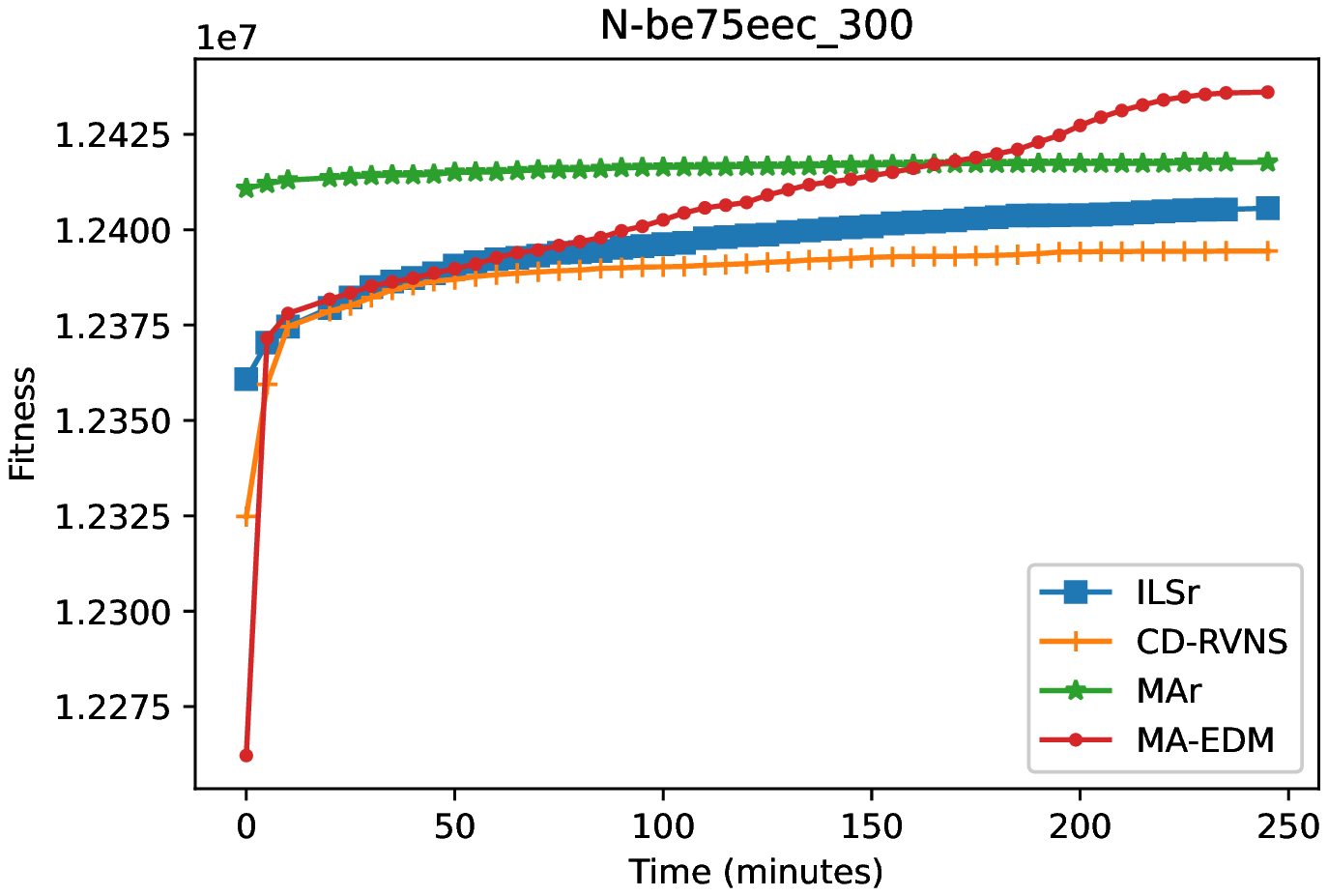}}
\subfloat[Nbe75eec500]{\label{fig:6}\includegraphics[width=0.48\textwidth]{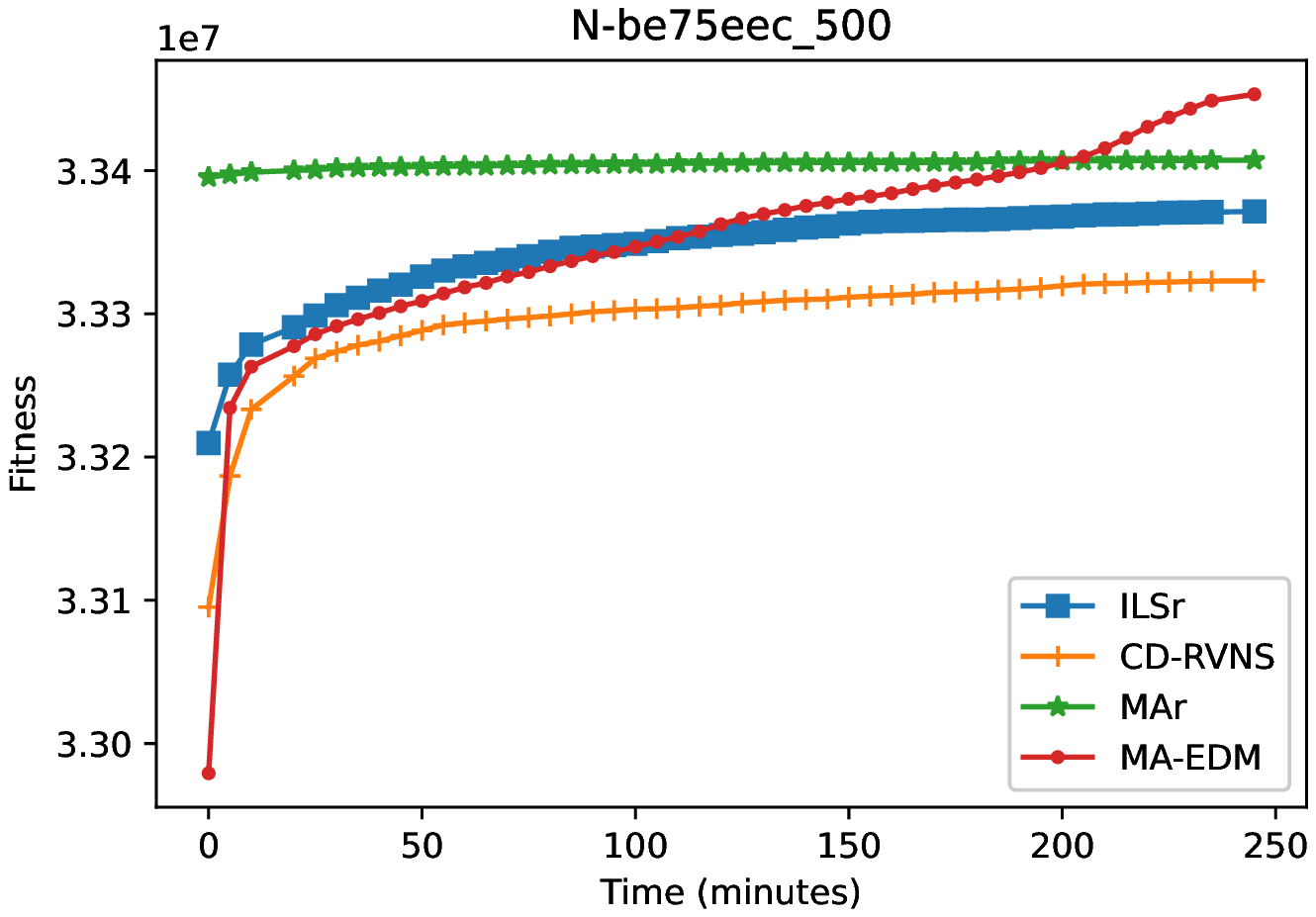}}
\qquad
\subfloat[Nbe75eec750]{\label{fig:7}\includegraphics[width=0.48\textwidth]{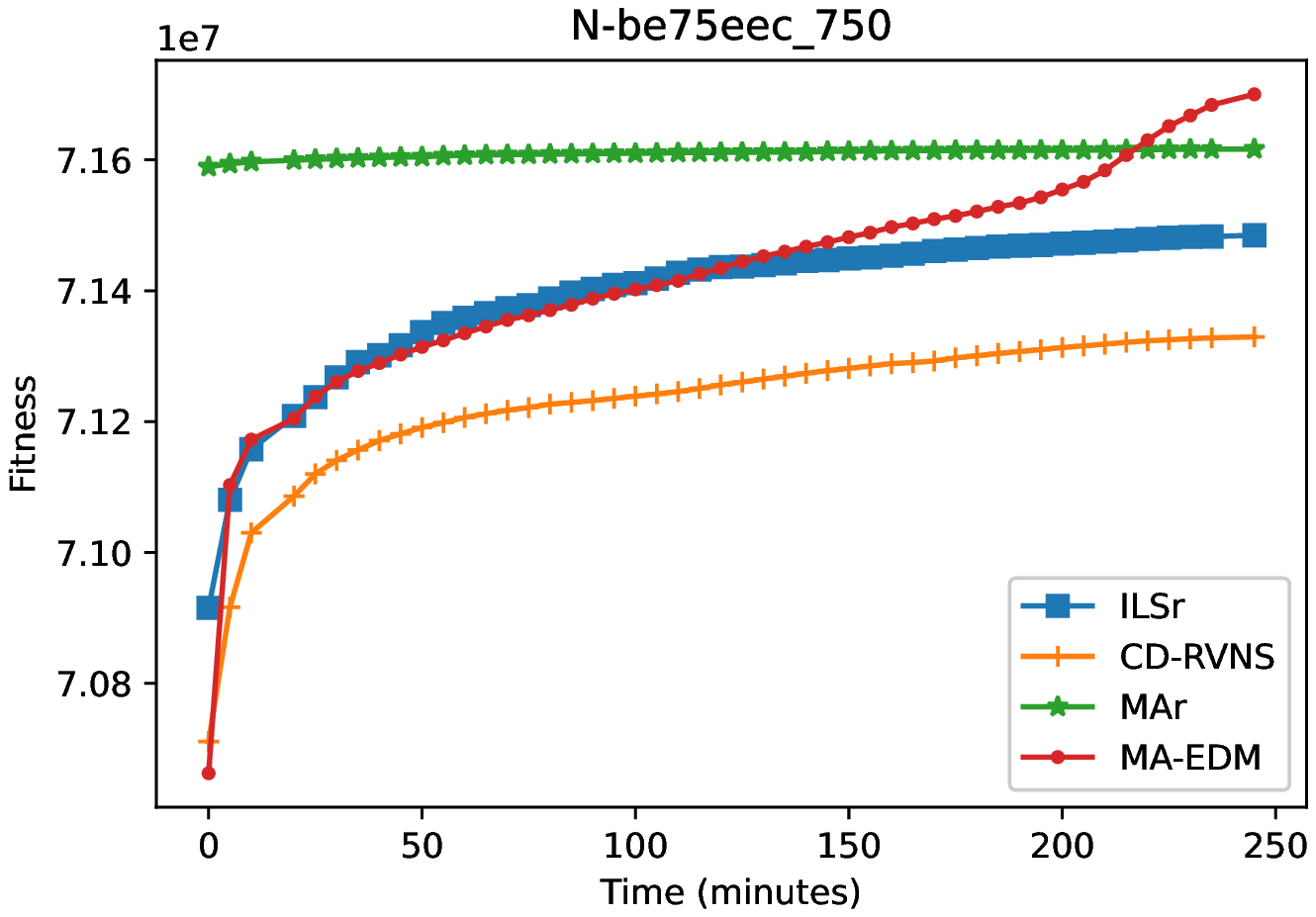}}
\subfloat[Nbe75eec1000]{\label{fig:8}\includegraphics[width=0.48\textwidth]{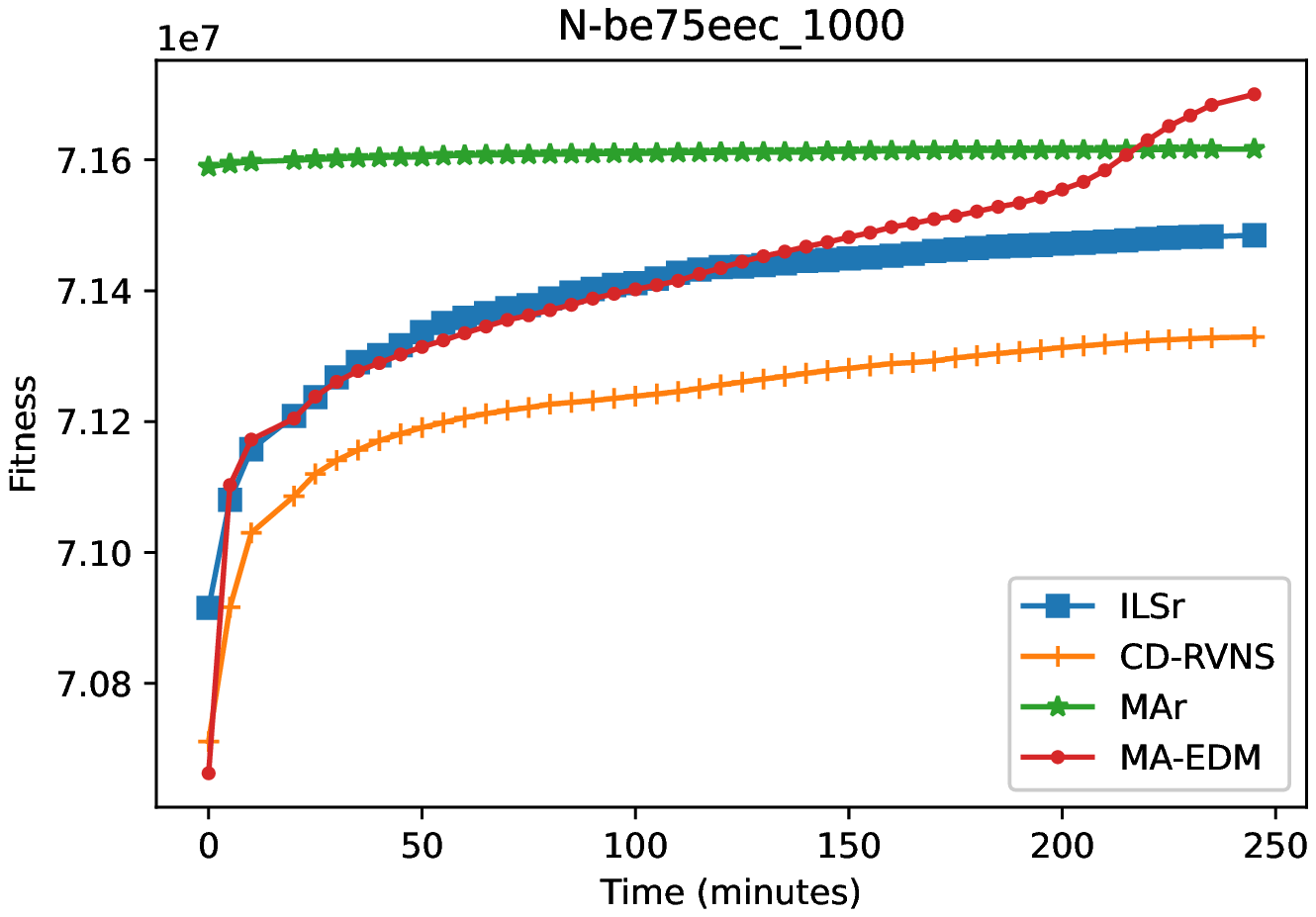}}
\qquad
\caption{Evolution of the best objective function attained during the run in four selected instances (mean of 30 executions)} 
\label{fig:fitness_evolution} 
\end{figure}

In order to complement the above findings, the trend in the best solution evaluated during the search is shown in Figure~\ref{fig:fitness_evolution}.
As in the previous experiment, the mean of 30 executions is shown and,
in this case, the trajectory-based methods are also included.
First, note that for these long-term executions, population-based approaches
are much more effective.
Second, it is also clear that MA-EDM exhibits a slow convergence, and specifically, when
compared to MAr, it shows a lower solution quality during most of the search.
However, as the end of the search nears and the diversity decreases, intensification is promoted and MA-EDM is able to profit from the promising regions
identified and yield much better final quality values.
Finally, it is also important to note that these results do not indicate that
MA-EDM requires more than 3 hours to improve on the results of MAr.
If a more restrictive computational budget were to be applied, MA-EDM would reduce
its diversity faster, so even when using shorter times, MA-EDM might also exceed MAr.

\section{Conclusions and Future Work}\label{sec5}

The LOP is a very popular NP-hard combinatorial optimization problem with numerous applications.
Metaheuristics have excelled in solving the LOP.
However, their performance degrades when considering large LOP instances.
Knowing near-optimal solutions is important to properly measure the performance of solvers.
This paper represents an effort to improve further the best-known solutions for standard benchmarks used in the validation of LOP solvers.
Currently, most state-of-the-art population-based metaheuristics designed for the LOP rely on simple schemes based on partial restarts to avoid premature convergence.
However, recent design trends in the field of MAs show that for difficult combinatorial optimization
problems, important advances can be achieved when more elaborate diversity management mechanisms are
included; 
specifically, those schemes that alter their internal components to relate the amount of diversity
maintained to the ratio between the elapsed execution period and the stopping criterion, have
yielded significant improvements.
In order to design a LOP solver that is able to improve further the new best-known solutions,
a novel MA that includes the mentioned design principle is devised.
This proposal is called the Memetic Algorithm with Explicit Diversity Management (MA-EDM).
MA-EDM incorporates a novel survivor selection strategy that is based on avoiding individuals 
that are too close to each other, which is done dynamically and considers the permutation deviation distance.

The experimental validation clearly confirms the advantages provided by this design principle for the LOP.
New best-known solutions were attained in 293 well-known instances, and the comparison against the state of the art
shows quite remarkable improvements.
Moreover, the advantages provided by MA-EDM were confirmed with proper statistical tests, and additional analyses were conducted
to understand the reasons behind the important advances in MA-EDM in comparison to alternative MAs.
The differences between the amount of diversity maintained in MA-EDM and other MAs is noticeable, which explains the important improvements contributed by MA-EDM. 
Thus, this research confirms the importance of proper diversity management in LOP solvers, and establishes a significant amount of new best-known solutions that should be taken into account in future validations of LOP solvers.

One of the most important drawbacks of our approach is that high-quality solutions are only attained near the end of the optimization process.
Introducing some alternative mechanisms that combine the intensification and exploration phases might improve the anytime behavior of MA-EDM~\cite{Radulescu:13}, which could allow its application in alternative conditions, where not just the quality of the final solution is taken into account.
Another topic that was not considered in this research was an analysis of the different ways of calculating distances between solutions.
Considering alternative metrics such as the exact match distance and the reversal distance seems quite promising.
Finally, the supplementary material associated with this research shows that as the dimensionality of the instances increases, the distance between the
best and mean objective function attained by MA-EDM grows.
Since the neighborhood considered in this paper exhibits quadratic growth, dealing with matrices with several thousand rows and columns might drastically limit the number of generations evolved, which could deteriorate the quality considerably.
In order to develop more robust and scalable approaches, techniques that have been devised specifically to improve
scalability might be incorporated.
Specifically, applying dimensionality reduction techniques and parallelization might allow the attainment of additional new best-known solutions.

\section*{Acknowledgments}
Authors acknowledge the financial support from CONACyT through the ``Ciencia B\'asica'' project no. 285599
and the support from ``Laboratorio de Super\-c\'om\-puto del Bajio'' through the project 300832 from CONACyT.

\bibliography{sn-bibliography}


\end{document}